\documentclass[12pt]{amsart}

\usepackage{a4wide}

\usepackage{amsmath,amscd}
\usepackage{amsfonts}
\usepackage{amsthm}
\usepackage{mathrsfs}
\usepackage{amssymb,amscd,latexsym}
\usepackage[title,titletoc,toc]{appendix}

\usepackage{pstricks,pst-node,pst-plot}
\usepackage{tikz}
\usetikzlibrary{calc}
\usetikzlibrary{arrows}

\usepackage{tikz}
\usetikzlibrary{positioning}

\tikzset{>=stealth}

\theoremstyle{plain}

\theoremstyle{definition}		% makes these non-italic

\usepackage{verbatim}

\newcommand{\Z}{\mathbb{Z}}
\newcommand{\N}{\mathbb{N}}
\newcommand{\Q}{\mathbb{Q}}
\newcommand{\C}{\mathbb{C}}
\newcommand{\R}{\mathbb{R}}
\newcommand{\F}{\mathbb{F}}
\newcommand{\pr}{\mathrm{pr}}
\newcommand{\Spec}{\mathrm{Sp}}
 \newcommand{\ii}{\mathbf{i}} 
  \newcommand{\cl}{\mathrm{cl}} 
\newcommand{\dn}{\mathrm{dn}}
\usepackage{pstricks,pst-node,pst-plot}
\usepackage{graphicx}
\title{Model theory and geometry of representations  of rings of integers}
\author{B.Zilber and L. Shaheen}
\email{}
\address{Mathematical Institute, University of Oxford, Oxford}
\date{21 September 2016}
\begin{document}
\maketitle

%\address{Mathematical Institute, University of Oxford, Oxford}

\begin{abstract}

The aim of this project is to attach a geometric structure to the  ring of integers. It is generally assumed that the   spectrum  $\mathrm{Spec}(\mathbb{Z})$ defined by Grothendieck serves  this purpose. However, it is still not clear what geometry this object carries.  A.Connes and C.Consani published recently an important paper  which introduces a much more complex structure called {\em the arithmetic site} which includes $\mathrm{Spec}(\mathbb{Z}).$

Our approach  is based on the generalisation of constructions applied by the first author for similar purposes in non-commutative (and commutative) algebraic geometry.

The current version is quite basic.  We describe a category of certain representations of integral extensions  of $\Z$ and establish its tight connection with the space of elementary theories of pseudo-finite fields. From model-theoretic point of view the category of representations is a
multisorted structure which we prove to be superstable  with pregeometry of trivial type. It comes as some surprise that a structure like this can code a rich mathematics of pseudo-finite fields.

\end{abstract}
\maketitle

\section{Introduction}\label{sec:intro}

The idea that the integers $\Z$, as an object of number-theoretic studies, could be better understood by associating a geometric object  to it is an accepted point of view. Following A.Grothendieck one views $\Z$ as a potential ``coordinate algebra'' of a geometric object which one can refer to as $\mathrm{Spec}(\Z).$ Conventionally, one thinks of   $\mathrm{Spec}(\Z)$ as a set of prime ideals of ring $\Z$ but it is not clear what {\em structure}, that is what  relations and operations this set naturally acquires. As Yu.Manin points out in  \cite{manin}, it is not understood  even  what     $\mathrm{Spec}(\Z)\times \mathrm{Spec}(\Z)$  
is.  Manin also speculates in  \cite{manin}
on what  $\dim \mathrm{Spec}(\Z)$ could be (pointing to  three possible
answers arising in discussions, $1,3$ and $\infty$). 

A.Connes and C.Consani in a series of papers that go back to the 1990's developed  a rich and interesting theory around this problem and more recently  (see \cite{CC}) introduced and studied a relevant structure which they called {\em the arithmetic site}.  Their construction builds the arithmetic site, following Grothendieck's prescriptions, as a topos in which ``points'' correspond to representations of the monoid $\N^\times$ of positive integers. The conventional   $ \mathrm{Spec}(\Z)$   can be seen embedded in the arithmetic site.

Our work makes use of the above ideas and visions but approaches the problem from a somewhat different systematic model-theoretic point of view. This approach has been presented in the general setting in \cite{Fun?} and in more concrete applications in \cite{QZG} for quantum algebras at roots of unity,  in \cite{QMcruz}, \cite{QM} for a quantum-mechanical setting and in ongoing projects in the setting of schemes and stacks as well as $C^*$-algebras. The leading idea of this approach is that the respective ``coordinate algebra'' encodes a syntactic description of a semantically presentable ``geometric structure'',  which is supposed to be a {\em geometry} in a sense  generalising the notion of Zariski geometry (see \cite{ZBook}) to the level of  corresponding formal languages. One obvious advantage of this model-theoretic approach is that the resulting object is a structure in the model-theoretic sense (usually on a multi-sorted universe) and so therein  direct products, projections, as well as various dimensions are well-defined.  Another advantage is the availability of the toolbox of model-theoretic methods (quantifier-elimination, stability theory and the theory of Zariski geometries).

\medskip

% Finally we have developed some comparisons between the model theoretic notions of  these representation of  monoids actions. 

\medskip

We start by associating with the ring of integers $O_K$ of a number field $K$ a very simple two-sorted structure on sorts $A_K$ and $\Spec_K,$ where the  standard universe for 
$\Spec_K$ is just the conventional $\mathrm{Spec}(O_K)$ as a set, and $A_K$ is a sort that is projected on $\Spec_K$ by a map $\pr,$ the fibres of the projection, $\pr^{-1}(\mathfrak{p})$ having the structure of a one-dimensional vector space over $\F_\mathfrak{p}= O_K/\mathfrak{p},$ which ignores the additive structure of the vector space, 
 ``monoid-representations of $O_K$''.

Objects $(A_K,\Spec_K),$ which we call {\em arithmetic planes} over respective number fields, are linked together into a category by
morphisms $$\pi_{K,L}: A_K\to A_L\mbox{ for }K\supseteq L,$$
which on each fibre $\pr^{-1}(\mathfrak{p})$ induce the norm map $\F_\mathfrak{p}\to \F_\mathfrak{q}$ for prime ideals
$\mathfrak{p}\supseteq \mathfrak{q}$ of respective rings.

This defines the multisorted structure, the object of our  model-theoretic analysis. Before stating results of this analysis we must note that it is not difficult to see that this structure is definable in the ring of finite adeles, well undersood object of model theory, in particular studied recently in much detail by J.Derakhshan and A.Macintyre, see \cite{DM} and their forthcoming papers.  These studies of the rings of adeles shed some light at our structure but do not explain the 
limits of expressive power of our language.  

Our main theorem states that the theory of the structure is superstable and allows elimination of quantifiers to certain family of core formulas of geometric flavour. This is not very surprising given that the structure is defined in terms of very simple relations, and indeed the only definable subsets on each spectral sort $\Spec_K$ are certain unary predicates. However, the language of the structure has quite a considerable expressive power:
we prove that to any point $v\in \Spec_K$ in a   model of the theory one can associate a pseudo-finite field $\F_v$ so that 
$$\mathrm{tp}(v)=\mathrm{tp}(w)\mbox{ if and only if }\F_v\equiv \F_w,\ \ v,w\in \Spec_K$$
where the elementary equivalence is in the language of fields extended by the names of its algebraic elements.

The quantifier-elimination theorem allows to define a natural topology on the multisorted structure and determine the dimensions (which we take to be just the $U$-ranks) of closed, and more generally definable, subsets. In particular,
$$\mathrm{U}(\Spec_K)=1,\ \ \mathrm{U}(A_K)=2,\ \mathrm{U}(\pr^{-1}(p))=1$$
for any number field $K$ and any $p\in \Spec_K.$ 
If we take the dimension to be the Morley rank then  the first two values are equal to $\infty$ and $\mathrm{MR}(\pr^{-1}(p))$ is finite (so our analysis explains two of the three versions of dimension of $\mathrm{Spec}(\Z)$ suggested by Yu.Manin in \cite{manin}). 

Another interesting result is that a large subalgebra of the boolean algebra of definable subsets on a sort $\Spec_K$ can be given a probabilistic measure, which is just the {\em natural} or {\em analytic  density} in the sense of number theory. We don't know if this measure is well-defined on all definable subsets.

We note that the topology is not Noetherian and that in the standard model the arithmetic planes $A_K$ and the spectral line $\Spec_K$ are not compact. However, we find that there are compact models and determine the minimal compact model. This model has finitely many non-standard (infinite) primes in each $\Spec_K,$ more precisely, the number is equal to $\deg K/\Q,$ and these primes $w$ are characterised by the property that any polynomial over $\Z$ splits into linear factors in $\F_w.$     

\medskip

{\bf Acknowledgement.} We would like to express our gratitude to J.A. Cruz Morales who read the paper and made many useful comments. 

\section{A bundle over the spectrum of $O_K$}\label{spec}

\bigskip

All fields $K$ below are number fields and rings are $O_K,$ the integers of the fields.  We generally assume that the fields belong to a collection  $\mathcal{R}$ closed under intersections such that for any $K,L\in   \mathcal{R}$ an embedding $L\subset K$ is Galois. By default we assume that the minimal object in $\mathcal{R}$ is $\Q$ but in general it could be  any number field.

 The spectrum of a ring $O_K,$ $\mathrm{Spec}(O_K),$
 is the collection of all prime ideals of $O_K$. Equivalently, for Dedekind rings, the collection of maximal ideals $\mathrm{maxSpec}(O_K),$
 together with the zero ideal. We  denote for brevity $$\Spec_K: =\mathrm{maxSpec}(O_K),  $$
 equivalently, the collections of irreducible representations of $O_K.$ This will be our main universe, a geometric space (of closed points), which agrees with the notion of universe for a Zariski geometry, see \cite{ZBook}. Note that when we introduce a ``Zariski'' topology on  $\Spec_K$ this universe as a whole is assumed to be closed, which in scheme-theoretic terms amounts to take the zero ideal into account.

 The maximal (that is non-zero prime) ideals we call points of $\Spec_K$ and often just say points of
  the spectrum.

 However, we want to  consider the  $O_K$  as monoids and to every point $p$ of $\Spec_K$ we put in correspondence an irreducible  $O_K$-module  $\mathrm{Fb}_{p}$
 which we consider as a monoid-module, that is we ignore the additive structure on $\mathrm{Fb}_{p}$ but do distinguish the zero element $0\in \mathrm{Fb}_{p}.$
  
  Since $K$ acts on $\mathrm{Fb}_{p}$ as $O_K/p$ and 
   $p$ is maximal,  $O_K/p\cong \mathbb{F}_q,$ a finite Galois field, $q=q(p)=\ell^n$ for some prime $\ell$ which depend on $K$ and $p.$ So, choosing a non-zero $a_p\in \mathrm{Fb}_{p}$    we can identify $\mathrm{Fb}_{p}$ with
     $\F_q.a_{p},$ where the action of $\mathbb{F}_q^*$ on $\mathrm{Fb}^{*}_{p}=\mathrm{Fb}_{p}\setminus \{ 0\}$ is free. 
     
     The characteristic property of $\mathrm{Fb}_{p}$ is that $$\mathrm{An}_{O_K}(\mathrm{Fb}_{p})=p,$$  
     the annihilator of $\mathrm{Fb}_{p}$ (in fact, of any non-zero point in $\mathrm{Fb}_{p}$) is equal to the ideal in $O_K$ generated by $p.$ 

\subsection{The arithmetic plane over $K$}     
     We define the 2-sorted (standard) universe  $\left( A_K, \Spec_K\right)$ with the projection
     $$\pr: A_K\setminus \{ 0\}\to \Spec_K$$ 
     where 
    $$A_K:= \bigcup_p \mathrm{Fb}_{p}\mbox{ and } \pr(x)=p \leftrightarrow x\in \mathrm{Fb}_{p}.$$ 
    
see the picture for $K=\Q,$ $O_K=\Z$ below.
% in \ref{1.2}.
\bigskip

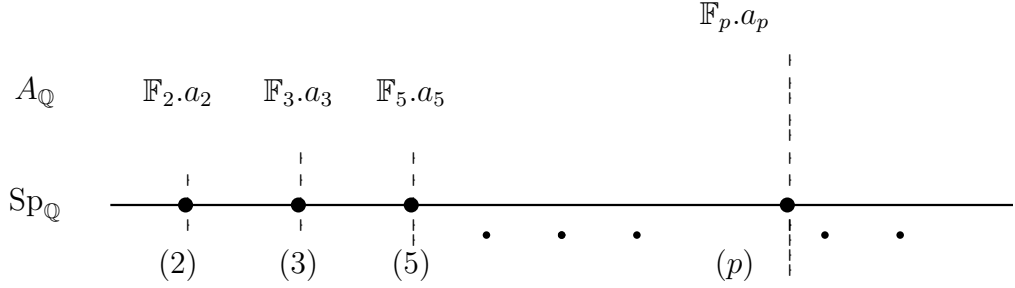
\begin{tikzpicture}
%\draw [line width =0.5mm] (1,0.6) node {$\mapstochar$};
\draw (1,0.3) node {$\mapstochar$};
\draw (1,0.0) node {$\mapstochar$};
\draw (1,-0.3) node {$\mapstochar$};
%\draw (1,-0.5) node {$\mapstochar$};
%%%%%%%%%%%%%%%%

\draw [line width =0.5mm] (2.5,0.6) node {$\mapstochar$};
\draw (2.5,0.3) node {$\mapstochar$};
\draw (2.5,0.0) node {$\mapstochar$};
\draw (2.5,-0.3) node {$\mapstochar$};
%\draw (2.5,-0.5) node {$\mapstochar$};

%%%%%%%%%%%%%%

\draw[line width =0.3mm] (0,0)--(2,0)--(5,0)--(7,0)--(10,0)--(11,0)--(12,0);

\fill (1,0) circle (0.1);
\fill (2.5,0) circle (0.1);
\fill (4,0) circle (0.1);

\fill (9,0) circle (0.1);

\draw (0.9,1.5) node {$\mathbb{F}_{2}.a_{2}$};
\draw (0.9,-0.8) node {$(2)$};

\draw(8.3,-0.8) node {$(p)$};
\draw(2.5,1.5) node {$\mathbb{F}_{3}. a_{3}$};
\draw (2.5,-0.8) node {$(3)$};
\draw (4,-0.8) node {$(5)$};

\draw (-1,0) node {$\Spec_{\Q}$};

%\draw [line width =0.5mm] (4,0.6) node {$\mapstochar$};
%\draw [line width =0.5mm] (4,0.6) node {$\mapstochar$};
\draw [line width =0.5mm] (4,0.6) node {$\mapstochar$};
\draw (4,0.3) node {$\mapstochar$};
\draw (4,0.0) node {$\mapstochar$};
\draw (4,-0.3) node {$\mapstochar$};
\draw (4,-0.5) node {$\mapstochar$};

\draw(4,1.5) node {$\mathbb{F}_{5}.a_{5}$};

\draw [line width =0.5mm] (9,0.6) node {$\mapstochar$};
\draw [line width =0.5mm] (9,.3) node {$\mapstochar$};
\draw [line width =0.5mm] (9,.9) node {$\mapstochar$};
\draw (9,-0.3) node {$\mapstochar$};
\draw (9,-0.3) node {$\mapstochar$};
\draw (9,-0.5) node {$\mapstochar$};
\draw (9,-0.7) node {$\mapstochar$};
\draw (9,-0.9) node {$\mapstochar$};

\draw (9,1.1) node {$\mapstochar$};

\draw (9,1.4) node {$\mapstochar$};

\draw(8.3,2.5) node {$\mathbb{F}_{p}.a_{p}$};

\draw (9,1.6) node {$\mapstochar$};

\draw (9,1.9) node {$\mapstochar$};

\fill (5,-0.4) circle (0.05);
\fill (6,-0.4) circle (0.05);

\fill (7,-0.4) circle (0.05);

\draw(-1,1.5) node {$A_{\Q}$};
\fill (9.5,-0.4) circle (0.05);
\fill (6,-0.4) circle (0.05);

\fill (10.5,-0.4) circle (0.05);

\end{tikzpicture}

\bigskip

\bigskip

We think about $A_K$ as an {\bf arithmetic plane over $K$} and $\Spec_K$ the horizontal axis of the plane.

\medskip

Thus,  the $\mathrm{Fb}^{*}_{p}$ are fibres of the projection. Note that for each standard fibre \\ $O_K\cdot a_p=\mathrm{Fb}_{p},$ but this is not  first-order expressible uniformly in $p.$  

We call the structure described above $\mathrm{Rep}(O_K),$ the representations of the monoid $O_K.$

We summarise: the language $L_K$ of $\mathrm{Rep}(O_K)$ has two sorts, $A_K$ and $\Spec_K,$ names for unary operations $m: x\mapsto mx$ for
elements $m\in K,$ and the map $\pr: A_K\to \Spec_K.$ 

\medskip

{\bf Remark.} We can equivalently represent $\mathrm{Rep}(O_K)$ in a one-sorted way, with just a sort
$A_K$ with equivalence relation $E$ instead of $\pr,$ such that the equivalence classes are exactly the $\mathrm{Fb}_{pK}.$

\subsection{Topology} \label{top}
To every ideal $m\subset O_K$ we associate the $\emptyset$-definable subset\\ $\mathcal{S}_{m}= \{ x\in A_K: \mathrm{An}(x)\supseteq m\}$, equal to  
 $\bigcup_{p|m} \mathrm{Fb}_{p}.$  Since every $m$ is finitely generated, the latter is a union of finitely many orbits. 
 
 We call such sets and finite union of those {\bf closed in $A_K.$}  
 
 Along with this topology on $A_K$ we define a topology on $\Spec_K$ with the basis of closed sets of the form $\pr(\mathcal{S}_{m}).$
 
 We will refer to as the {\em conventional} Zariski topology on $A_K$ and $\Spec_K.$

It is not hard to prove the following statement which will be superseded by Theorem~\ref{ThmQE}.

\medskip

{\bf Proposition.} {\em The complete first-order theory of $\mathrm{Rep}(O_K)$ has an explicit axiomatisation $T(K).$ The theory admits elimination of quantifiers.
This theory is $\omega$-stable of finite Morley rank. The sort $\Spec_K$ is strongly minimal  and of trivial type. The sort $A_K$ is of Morley rank $2$ and the fibres of $\pr$ are either finite or strongly minimal.}

\section{The Multi-Sorted Structure of Representations}\label{SR}

\subsection{The language} \label{LR}
In this section we present a construction of  a multi-sorted structure with sorts of the form $A_K, \Spec_K$ where $K$ runs through a family $\mathcal{R}$ of (by default all) number fields. 

The language $L_\mathcal{R}$ of the {\bf multisorted representation structure} $\mathrm{Rep}(\mathcal{R})$ 
will be the union of the languages $L_{K},$ $K \in \mathcal{R},$
 extended by maps $\pi_{K,L}: A_K\to A_L,$ for $K,L \in \mathcal{R},$
$\mbox{O}_{L}\subseteq \mbox{O}_{K}.$

We will also use notation for definable maps $$\pi^\Spec_{K,L}: \Spec_{K}\to \Spec_{L}$$
defined by  $$\pi^\Spec_{K,L}(p{K})=p{L} \Leftrightarrow \mathrm{Fb}_{pK}\subseteq \pi_{K,L}(\mathrm{Fb}_{pL})$$
\medskip

\subsection{Notation} \label{pi}
%We will now describe the maps $\pi_{\alpha,\beta}.$

\medskip

Let $\mathfrak{q}\subset O_K$ and  $\mathfrak{p}\subset O_L$
are prime ideals of the corresponding rings $O_L\subseteq O_K$ such that 
 $\mathfrak{q}\supset  \mathfrak{p}.$
In this case we usually use notations  
$p{K}:= \mathfrak{q}$ and  $p{L}:= \mathfrak{p}$ 
thus we have by this notation 
\begin{equation} \label{eq1}
p{K}\supseteq p{L},
\end{equation}
that is always
``$p{K}$ lies over $p{L}$ in $\mbox{O}_{K}$''. 
%In particular, if $p\in \Spec_{\Z}$ is a prime ideal in $\Z$ then $p{L}$ lies over $p.$ 

Recall that the ring homomorphism 
$$\mathrm{res}_{p{K}}:\mbox{O}_{K}\to \mbox{O}_{K}/p{K}=: \F_{p{K}}$$
can be seen as a residue map (or place) for a valuation of $K$ with the value ring equal to $\mbox{O}_{K}$ and the valuation ideal equal to $p{K}.$ Here  $ \F_{p{K}}\cong\F_{p^m}$ for some positive integer $m,$ which depends on $K.$ We define
$$|p{K}|:=p^m.$$

Note that if $K:L$ is Galois, acting by a Galois automorphism on  $K$ we move $p{K}$ to some $p'{K}\subset \mbox{O}_K.$  In this sense $pK$ runs through all prime ideals of $O_{K}$ which are Galois-conjugated to  $pK,$ while   $ \F_{p'{K}}\cong\F_{p^m}$ with the same $m.$

%\mbox{ and } \mathrm{res}_{p_\beta}:\Z[\beta]\to \Z[\beta]/(p_\beta)\cong \F_{q_{\beta}}$$ for some finite fields $ \F_{q_{\alpha}}$ and $\F_{q_{\beta}}. $ Under our assumptions 
%$ \F_{q_{\alpha}}\supseteq \F_{q_{\beta}}. $

We will often refer to $\mathrm{res}_{p{K}}$ as a {\bf naming homomorphism} as it associates elements $\gamma\in  \mbox{O}_{K}$ (which we see as ``names'') with elements in $\F_{p{K}}.$

\subsection{Orbits} \label{orbits} For each prime ideal $p{K}$ we associate a unique orbit $\mathrm{Fb}_{pK}$ on which
$O_{K}$ acts, for each $\gamma\in O_{K}\setminus p{K}$ and $x\in \mathrm{Fb}_{pK}$
$$ x\mapsto \gamma\cdot x\in \mathrm{Fb}_{pK}$$
and $\gamma\cdot x=0$ (a common zero) iff $\gamma\in p{K}.$ 
By definition we assume 
$$\gamma_1\cdot x=\gamma_2\cdot x\mbox{ iff } \gamma_1-\gamma_2\in p{K}$$
and we assume that the action is transitive. Thus, by definition, for any $\gamma\in O_{K}$ and $a\in \mathrm{Fb}_{pK},$

\begin{equation} \gamma\cdot a=\bar{\gamma}_{p{K}}\cdot a
\end{equation}
where $\bar{\gamma}_{p{K}}\in \F_{p{K}},$
the element with the name $\gamma.$ 
\medskip

\subsection{Corollary} For $q,q'\in \Spec_{K}$
$$q=q'\mbox{ {\bf iff} for any }\gamma_1,\gamma_2\in O_{K}\, \forall x\in \mathrm{Fb}_{q^K}\, \forall x'\in \mathrm{Fb}_{q'^K}\ \gamma_1\cdot x=\gamma_2\cdot x \leftrightarrow  \gamma_1\cdot x'=\gamma_2\cdot x'.$$

Note also that, given $a\in  \mathrm{Fb}_{pK}$ we can represent any element $b\in \mathrm{Fb}_{pK}$ as
$$b=\bar{\gamma}_{p{K}}\cdot a, \mbox{ for some } \bar{\gamma}_{p{K}}\in \F_{p{K}}.$$

We will also write accordingly
$$\mathrm{Fb}_{pK}= \F_{p{K}}\cdot a_{pK},$$ where $pK\mapsto a_{pK}$ is a cross-section $\Spec_{K} \to A_{K}$

\label{morphisms}\subsection{Morphisms between sorts}
Let $\mathrm{Norm} _{K,L}:
\F_{p{K}}\to \F_{p{L}} $ be the norm map.  Define,
the morphisms
$$\pi_{K,L}: A_{K}\to A_{L}$$
by its action on each orbit
$$\pi_{K,L}:\mathrm{Fb}_{pK}\to \mathrm{Fb}_{pL},$$
for each prime $p{K}$ over $p{L},$
which is defined once the sections $a_{pK}\in \mathrm{Fb}_{pK}$ and $a_{pL}\in \mathrm{Fb}_{pL}$ are given:
 for $\eta\in  \F_{p{K}}$
$$\pi_{K,L}(\eta\cdot a_{pK}):= \mathrm{Norm} _{K,L}(\eta)\cdot a_{pL}.$$

Note that this definition applies to  all primes $p'K$ in $O_K$ over $pL$ and morphisms
 $$\pi_{K,L}:\mathrm{Fb}_{p'{K}}\to \mathrm{Fb}_{pL}.$$

\subsection{Remark} It is crucial for the definition of  morphisms $\pi_{K,L}$ that we treat the $O_K$ as monoids (not rings), that is there is no additive structure on the fibres $\mathrm{Fb}_{p{K}}.$

 \bigskip

{\bf{\underline{Naming homomorphisms}}}:

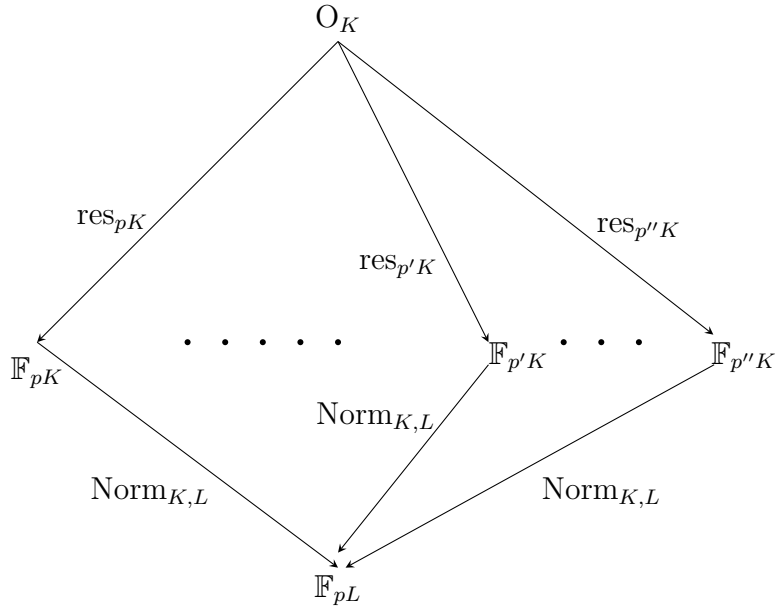
\begin{center}

\begin{tikzpicture}
\draw[-] (4,0.3) node {$\mbox{O}_{K}$};

\draw[-] (0,-4.4) node {$\F_{p{K}}$};
\draw[-] (6.4,-4.2) node {$\mathbb{F}_{p'{K}}$};
\draw[-] (9.4,-4.2) node {$\mathbb{F}_{p''{K}}$};
\draw[-] (1,-2.4) node {$\mbox{res}_{p{K}}$};
\draw[-] (4.8,-3) node {$\mbox{res}_{p'{K}}$};
\draw[-] (8,-2.5) node {$\mbox{res}_{p''{K}}$};
\draw[-] (4,-7.3) node {$\F_{p{L}}$};
\draw[-] (1.5,-6) node {$\mbox{Norm}_{K,L}$};
\draw[-] (4.5,-5) node {$\mbox{Norm}_{K,L}$};
\draw[-] (7.5,-6) node {$\mbox{Norm}_{K,L}$};

\draw[->] (4,0) -- (0,-4);
\draw[->] (0,-4) -- (4,-7);

\draw[->] (4,0) -- (8.99,-3.9);

\draw[->] (4,0) -- (6,-4);

\draw[->] (6,-4.3) -- (4,-6.8);

\draw[->] (9,-4.3) -- (4.1,-7);

\fill (2,-4) circle (0.04);

\fill (2.5,-4) circle (0.04);

\fill (3,-4) circle (0.04);

\fill (3.5,-4) circle (0.04);

\fill (4,-4) circle (0.04);

\fill (7,-4) circle (0.04);

\fill (8,-4) circle (0.04);

\fill (7.5,-4) circle (0.04);

\end{tikzpicture}

\end{center}

%\includegraphics[width=20cm, height=10cm]{naming.pdf}%
%\medskip
\subsection{Lemma}   
%Since  $\mathrm{Norm}_{\alpha,\beta}$
  %is a surjective homomorphism of monoids identical on  $\F_{q_\beta}$, the map 
{\em  The action of $\mathrm{O}_{K}$ on $\mathrm{Fb}_{pK}$ induces via
$\pi_{K,L}$ an action of $\mathrm{O}_{K}$ on $\mathrm{Fb}_{pL}.$ 

This action extends the action  of $\mathrm{O}_{L}$ on $\mathrm{Fb}_{pL}.$

This action does not depend on the sections $a_{pK}$ and $a_{pL}.$}

{\bf Proof.}  Given $\gamma\in \mathrm{O}_{K}$ and $b\in \mathrm{Fb}_{pL}$
 define
 \begin{equation} \label{actiona}
 \gamma\cdot_{pK} b:=\mathrm{Norm} _{K,L}(\bar{\gamma}_{pK})\cdot b.
 \end{equation}
 
This is clearly an action since norm is multiplicative. We will write it without the subscript $p{K}$ when the latter is clear from the context.

For $\gamma\in \mathrm{O}_{L}$ we have $\bar{\gamma}_{pK}\in \F_{pL}$ and $\bar{\gamma}_{pK}=\bar{\gamma}_{pL}$ since $pK\cap \mathrm{O}_{L}=pL.$ Hence $\mathrm{Norm} _{K,L}(\bar{\gamma}_{pK})=\bar{\gamma}_{pL}$ for such $\gamma$ and so this is its usual action. Clearly, it does not depend on the choice of sections.

Now we show that the $\pi_{K,L}$ preserve the action by 
$\mathrm{O}_{K},$ that is 
\begin{equation}
\pi_{K,L}(\gamma\cdot a_{pK})=
\gamma\cdot\pi_{K,L}(a_{pK})=\gamma\cdot a_{pL}.
\end{equation}  

Indeed, by definition $\gamma\cdot a_{pL}=\mathrm{Norm} _{K,L}(\bar{\gamma}_{pK})\cdot a_{pL},$ $\pi_{K,L}(\bar{\gamma}_{pK}\cdot a_{pK})=\mathrm{Norm} _{K,L}(\bar{\gamma}_{p{K}})\cdot a_{p{L}} $
and $\gamma\cdot a_{p{K}}= \bar{\gamma}_{pK}\cdot a_{pK}.$ The equality follows.\qed

\medskip

\subsection{ Remark} The definition of the action can equivalently be written as 
$$\gamma\cdot_{p{K}} b:=\mathrm{Norm} _{K,L}(\mathrm{res}_{pK}(\gamma))\cdot b.$$

\subsection{Corollary} {\em The actions  $\cdot_{pK}$ and   $\cdot_{p'{K}}$ on $\mathrm{Fb}_{pL}$ 
by $\mbox{O}_{K}$  coincide if and only if $p'{K}=p{K}$} 

Indeed, to see that the actions differ when  $p'{K}\neq p^K$ consider a $\gamma\in  p'{K}\setminus pK.$ Then by definition for any $b\in \mathrm{Fb}_{pL},$
$\gamma\cdot_{p'{K}}b=0$ and $\gamma\cdot_{p{K}}b\neq 0.$

\subsection{Remark} Note also that  $\mathrm{Norm}_{K,L}$ can be alternatively defined as the map
$$\mathrm{Norm}_{K,L}(x)=   x^{\frac{|p{K}|-1}{|p{L}|-1}}.$$

It follows that for a non-zero $y\in \mathrm{Fb}_{pL}$
\begin{equation} \label{eq2}
|\pi_{K,L}^{-1}(y)\cap \mathrm{Fb}_{pK}|=\frac{|pK|-1}{|pL|-1}.
\end{equation}

Now note that one can determine, for $y\in A_{L},$
 the number  $|\pi_{K,L}^{-1}(y)|$ once one knows  that
  $y\in \mathrm{Fb}_{pL}$ and one knows all the $p'K$ which lie over $pL$ in $\mathrm{O}_{K}.$

\subsection{Special predicates on $\Spec_{L}$}
We consider the $\pi_{K,L}$-{\bf multiplicities} of primes $pL\in \Spec_L.$
%points $y\in A_{L}$ to be the numbers $  |\pi^{-1}_{K,L}(y)|.$
%For positive integers $N$ we set 
% $$\psi_{K,L}^N:=\{ y\in A_{L}:    |\pi^{-1}_{K,L}(y)|\le N\},$$ definable subsets of $A_{L}.$ 

Define for positive integers $N$ the binary relation $P^N_{K,L}\subset  \mathrm{O}_{K}\times \mathrm{O}_{L}$ by the first order formula:

$$P^N_{K,L}(pK,pL)\equiv  \exists y\in \mathrm{Fb}_{pL}\  |\pi_{K,L}^{-1}(y)\cap \mathrm{Fb}_{pK}|\le N.$$

Note that by the definition of $\pi_{K,L}$ it follows that $P^N_{K,L}(pK,pL)$ implies  $pL\subseteq pK.$

%$$\Psi_{K,L}^N:=%\pr(\psi_{K,L}^N)=
%\{ p\in \Spec_{L}:\exists y\in A_{L} \left( |\pi^{-1}_{K,L}(y)|\le N \ \& \ \pr(y)=p\right)\}$$

%We definably extend the language $L_\mathcal{R}$ to $L_\mathcal{R}^\sharp$
% by adding to it unary  predicates $\Psi_{\alpha,\beta}^N$ and $\psi_{\alpha,\beta}^N$    on the sorts $\Spec(\beta)$ and $A(\beta)$ respectively,

%$\Psi_{\alpha,\beta}^N$ and  $\Phi_{\alpha,\beta}^d$
%for all positive integers $N.$ 
%and $d$ definable in $L_\mathcal{R}$ as follows

%$$\Psi_{\alpha,\beta}^N(p)\equiv \exists y\in A(\beta) \left( |\pi^{-1}_{\alpha,\beta}(y)|\le N \ \& \ \pr(y)=(p)\right)$$%\equiv \forall y\in  A(\beta) \left(\pr(y)=(p)\to   |\pi^{-1}_{\alpha,\beta}(y)|\le N\right).$$

%$$\Phi_{\alpha,\beta}^d(p)\equiv |\{ r\in \Spec\, \Z[\alpha]: O_r\subseteq  \pi^{-1}_{\alpha,\beta}(O_p)\}|\ge d.$$

%Note that given $p\in \Z[\beta],$
%$O_r\subseteq  \pi^{-1}_{\alpha,\beta}(O_p)$ if and only if the ideal generated by $p$ in $\Z[\alpha]$  contains the prime ideal generated by $r$ in the same ring. Hence $d_p,$ the maximal $d$ for which  $\Phi_{\alpha,\beta}^d(p)$ can hold is less or equal to $\deg (\Z[\alpha]:\Z[\beta]).$

%Given $pL\in \Spec_{L}$ set
%$$N_{pL,K}=\min \{ N: pL\in  \Psi_{K,L}^N\}.$$ 

\subsection{Lemma}\label{3.11} {\em 
Let $L\subseteq K,$ $pL\in \Spec_{L},$ $ pK\in \Spec_{K},$ and assume that $P^N_{K,L}(pK,pL)$ holds.
 Then
 
$\mathrm{(i)}\  N\ge \frac{|pL|^{|\F_{pK}:\F_{p{L}}|}-1}{|pL|-1}.$

{\rm (ii)}  In particular, there exists $N$ such that  $P^N_{K,L}(pK,pL)$ holds  
 for infinitely many  $pL\in \Spec_{L}$ with some extension  $pL\subseteq pK\in  \Spec_{K}$ 
 if and only if $|\F_{pK}:\F_{p{L}}|=1$ for all but finitely many of pairs $pK,pL.$ 
 
{\rm (iii)}  If there exists $N$ such that  $P^N_{K,L}(pK,pL)$ holds 
 for infinitely many  $pL\in \Spec_{L}$ with some extension  $pL\subseteq pK\in  \Spec_{K}$ 
 then $P^1_{K,L}(pK,pL)$ holds for all but finitely many pairs $pK,pL$ satisfying $P^N_{K,L}(pK,pL).$
}
\medskip

{\bf Proof.} (i) This quantative estimate is a direct consequence of
 (\ref{eq2}). 
 
(ii)  $|pL|$ is unbounded when $pL$ runs in an infinite subset of $\Spec_L.$

(iii)  $|\F_{pK}:\F_{p{L}}|=1$ is equivalent to the statement that $\pi_{K,L}$ induces a bijection $\mathrm{Fb}_{pK}\to \mathrm{Fb}_{pL}.$
 \qed
 %and simple facts about splitting polynomials in Galois extensions.

\subsection{Remark} \label{psi} The case $|\F_{pK}:\F_{p{L}}|=1$ corresponds to the fact that 
$\F_{pK}\cong\F_{pL}$ for the prime $pK\in \Spec_K$ over the  prime $pL\in \Spec_L.$ In case when $K$ is Galois over $L$ this means that 
the minimal polynomial of $\alpha$ over $L$ splits into linear factors modulo $pL.$ Then,  except for finitely many such $pL$ (over  which $pK$ ramifies) 
 there are exactly 
$ |K:L|$ distinct prime ideals $pK$
over $pL.$%   and thus $d_p=|K:L|.$ 

\subsection{}\label{3.14+}
Define  $$\Pi_{K,L}:=\{ pL\in \Spec_L: \ \exists  pK\in \Spec_{K}\ pL\subseteq pK\ \& \  \F_{pK}\cong\F_{pL}\}$$
$$\Psi_{K,L}:=\{ pL\in \Spec_L: \ \forall  pK\in \Spec_{K}\ pL\subseteq pK\ \to \  \F_{pK}\cong\F_{pL}\}$$
 
By \ref{3.11} both are first order definable (using $P^1_{K,L}(pK,pL)$).

\medskip

Note that  $pL\in \Psi_{K,L}$ if and only if every extension $pK$ of $pL$ to $\mathrm{O}_K$ splits completely. 

\medskip

We will show later (see  \ref{new} and    \ref{rem4.12})
that $\Psi_{K,L}$ can be expressed  essentially  in terms of $\Pi_{K,L}.$

\subsection{Lemma}\label{3.14++} {\em Suppose $K$ and $K'$ are  Galois conjugated extensions of $L.$ Then $$\Pi_{K,L}=\Pi_{K',L}.$$  
} 
\begin{proof} Let $g: K\to K'$ be an isomorphism over $L.$ Suppose $pL\in \Pi_{K,L}$ and so for some
$pK\in \Spec_K$ over $pL,$ 
$\mathrm{O}_K/pK\cong \mathrm{O}_L/pL.$ Then  $\mathrm{O}_{K'}/pK'=g(\mathrm{O}_K)/g(pK)\cong \mathrm{O}_L/pL.$ Hence $pL\in \Pi_{K',L}.$
\end{proof}

\subsection{Galois action on the spectra} Let $\alpha_1,\ldots,\alpha_k$ be the generators of the ring $O_K$ over $O_L$
 and
$\sigma\in \mathrm{Gal}(K:L).$ Let
$\alpha^\sigma_1,\ldots,\alpha^\sigma_k$ be the result of application of $\sigma$ to the generators.

Consider the formula with free variables $q, q'$ in sort $\Spec_K:$
\begin{equation} \label{qq}
\exists x,x'\in A_K:\ \pr(x)=q\ \&\ \pr(x')=q'\ \&\ 
\pi_{K,L}(x)=\pi_{K,L}(x')\ \&\ \bigwedge_{l=1}^{k}
\pi_{K,L}(\alpha_l x)=\pi_{K,L}(\alpha^\sigma_l x')
\end{equation}

{\bf Lemma.}\footnote{We should be able to write down a formula which defines an action of $\sigma$ on 
$A_K$ (not just $\Spec_K$) without the restriction to $\psi_{K,L}.$
}
 {\em Given
  $q\in \Spec_K$ such that
$\pi^\Spec_{K,L}(q)\in \Psi_{K,L}$ the formula (\ref{qq}) holds for $q'\in \Spec_K$
 if and only if $q'=q^\sigma.$}

{\bf Proof.} %We may assume $x=a_q$ and $x'=a_{q'}.$ 
Suppose (\ref{qq}) does hold. Then $p:=\pi^\Spec_{K,L}(q)=\pi^\Spec_{K,L}(q').$\footnote{In particular, we can write the formula (\ref{qq}) so that all the variables are relativised to $\cl(p).$}

Under our assumptions $\mathrm{Norm}_{\F_{q},\F_{p}}:\F_q\to \F_p$ is an isomorphism 
and 
 definition (\ref{actiona}) implies that for any $\alpha\in O_K,$ 
 $$\pi_{K,L}(\alpha\cdot x)=\bar{\alpha}_q\cdot y,\ \mbox{ for }y:=\pi_{K,L}(x),$$ 
 where $\bar{\alpha}_q=
\mathrm{res}_q(\alpha)\in\F_p.$

Since
$\pi_{K,L}(\alpha x)=\pi_{K,L}(\alpha^\sigma x'),$
we get that 
$$\bar{\alpha}_q=\bar{\alpha^\sigma}_{q'}$$
for all $\alpha\in O_K.$
In particular, $$\alpha\in \mathrm{\ker}(\mathrm{res}_q)
\Leftrightarrow
\alpha^\sigma\in \mathrm{\ker}(\mathrm{res}_{q'}).$$
But $\mathrm{\ker}(\mathrm{res}_q)=q$ and $\mathrm{\ker}(\mathrm{res}_{q'})=q',$ that is $q'=q^\sigma.$ \qed

\section{The quantifier elimination theorem}\label{QE}

\subsection{The space of finite fields} The multisorted structure $\mathcal{PF},$ {\em a space of finite fields,} consists of sorts $\mathcal{PF}_K,$   $\mathcal{AF}_K$ for number fields $K,$
with a surjective map $$\pr:    \mathcal{AF}_K   \to \mathcal{PF}_K.$$  

As was noted in subsection~\ref{pi}
 a  unital homomorphism 
 $\mathrm{res}:\mathrm{O}_K\to \F_{pK}$ can be seen as assigning names  in $ \mathrm{O}_K$ to elements in $\F_{pK}.$
 We will consider elements of   $ \mathrm{O}_K$ as extra constant symbols and
  $\F$  a structure in
 the  {\bf language of rings with names in $ \mathrm{O}_K.$} We denote the language $L_\mathrm{rings}(K).$
 
The fibre $\pr^{-1}(q)$
over a point $q\in  \mathcal{PF}_K,$  
has a structure given by the language $L_\mathrm{rings}(K)$ that can be identified 
with a finite field $\F_{q}$ (later, in \ref{UfD}, a pseudofinite field)
with names in $\mathrm{O}_K.$

We assume that for finite $\F_{q},$  
\begin{equation} \label{axiomFqq}
\F_{q}\cong \F_{q'}\Leftrightarrow q=q'
\end{equation}

We also assume that,

\begin{equation}\label{Fq}
\mbox{for every prime  }pK\subset \mathrm{O}_K\mbox{ there is }q\in  \mathcal{PF}_K\mbox{ such that }
\F_q=\F_{pK}:= \mathrm{O}_K/pK.
\end{equation}

There are maps $$\mathrm{Nm}_{K/L}:  \mathcal{AF}_K  \to  \mathcal{AF}_L,\ \ j_{K/L}:  \mathcal{AF}_L  \to  \mathcal{AF}_K\mbox{ for } \mathrm{O}_K\supseteq \mathrm{O}_L$$
between sorts, defined fibrewise as \begin{equation} \label{Nm}
  \mathrm{Norm}_{K/L}: \F_{{pK}}\to \F_{{pL}},\ \ j_{K/L}: \F_{{pL}}\to \F_{{pK}} \mbox{ for }   \F_{{pK}}\supseteq \F_{{pL}}.
  \end{equation}
 where $j_{K/L}$ is the canonical embedding.
 
 \medskip 
  Note that  \begin{equation} \label{Nmsurj}
   \mathrm{Nm}_{K/L}\mbox{ is surjective}
   \end{equation}
   since every prime $pL\subset \mathrm{O}_L$ can be lifted to a prime $pK\subset \mathrm{O}_K.$   
   
We also assume that   \begin{equation} \label{pLK}
\F_{{pK}}=\F_{pL}(\bar{\alpha}),\mbox{ for } \alpha \mbox{ such that }K=L(\alpha) \end{equation}
where $\bar{\alpha}=\mathrm{res}_{pK}(\alpha).$

\subsection{Remarks} 
 
(i)  Note that $\F_{pK}$ is determined by (\ref{pLK}) uniquely,  up to the naming, and we have 
one-to-one correspondence between extensions of $\F_{pL}$ in the language    $L_\mathrm{rings}(K)$ and the residue maps $\mathrm{res}_{pK}.$

(ii) The assumptions (\ref{axiomFqq})--(\ref{pLK}) are first-order axiom schemes of spaces of pseudofinite fields.  The axioms    (\ref{axiomFqq}) and  (\ref{Fq}) are written for each finite $q$ separately.   

 \subsection{Lemma} \label{srpf} {\em
There exists  a bijection $$\mathbf{i}:\  {A_K}\to \mathcal{AF}_K,\ \ 
{\Spec_K}\to  \mathcal{PF}_K,\ 
\mbox{ for all }K$$
between the (standard)  multisorted representation structure
$(A_K, \Spec_K)$  (section \ref{SR}) and the  space of finite pseudofinite fields  $(\mathcal{AF}_K, \mathcal{PF}_K)$  
such that  \begin{equation} \label{ii1}\ii\circ\pr=\pr\circ \ii,\ \ \mathrm{Nm}_{K/L}\circ \ii=\ii\circ \pi_{K,L}
 \end{equation} 
and for any $\gamma\in O_K, \ x\in {A_K}$ 
 \begin{equation} \label{ii2}
\ii(\gamma\cdot x)=\mathrm{res}_{\pr(x)}(\gamma)\cdot\ii(x). \end{equation} 

%If $\ii'$ is another such map then 
}

{\bf Proof.}  For each $L\subset K$ choose the sections $Sp_L\to A_L,$ $pL\mapsto a_{pL}.$ Then
$$x=\bar{\gamma}\cdot a_{pL}\Rightarrow \ii(x):=\bar{\gamma},\ \mbox{ for }x\in A_L,\ \gamma\in O_K,
\bar{\gamma}:=\mathrm{res}_{pL}(\gamma)\in \F_{pL}$$
determines the map with required properties. 
\qed
 
\bigskip

\subsection{Corollary (The space of pseudofinite fields).} \label{UfD} {\em  Given an  ultrapower $(^*A_K, {^*\Spec_K})$ (over an ultrafilter $\mathcal{D}$) of   the standard $(A_K, {\Spec_K})$ there  
exists a   space of pseudofinite fields  $({^*\mathcal{AF}}_K, {^*\mathcal{PF}}_K)$ and a bijection  $$\mathbf{i}:\  {^*A_K}\to {^*\mathcal{AF}}_K,\ \ 
{^*\Spec_K}\to  {^*\mathcal{PF}}_K,\ 
\mbox{ all }K$$ 
with the same properties (\ref{ii1}) and (\ref{ii2}) as above. More precisely, $({^*\mathcal{AF}}_K, {^*\mathcal{PF}}_K)$ can be constructed as the ultrapower of $({\mathcal{AF}}_K, {\mathcal{PF}}_K)$ over the same ultrafilter.}
  
\subsection{Algebraic closure}\label{cl}
We call a substructure $M$ of the multisorted structure $({\mathcal{AF}'}_K, {\mathcal{PF}'}_K)$ {\bf algebraically closed} if 

(i) all the finite primes and the fibres over them belong to $M;$

(ii) for any $pL\in \Spec_L\cap M$  we have $\pr^{-1}(pL)\subset M,$

(iii) and  for any $K\supseteq L$ every  $pK\in \Spec_K$ such that $pL\subseteq pK$ belongs to $M.$ 

\medskip

The same definition works for a substructure of the multisorted structure $({^*A}_K, {^*\Spec}_K).$

\medskip

Given arbitrary $X$ in a multisorted structure, we will write
$\cl(X)$ for the minimal algebraically closed substructure containing $X.$ 

\bigskip

{\bf Example.}  The algebraic closure of a point $x\in A_K$
is equal to the closure of the prime $pK\in \Spec_K,$
 such that $pK=\pr(x).$ In its turn $\cl(pK)$ contains, for every 
$L\subseteq K,$ the point $pL=\pi^\Spec_{K,L}(pK).$ But, at the same time, by definition $pK\in \cl(pL).$ Summarising,
one gets $\cl(x)=\cl(pK)$ by starting from $p\in \Spec_\Q,$ the unique prime in $\Z$ under $pK,$ and then adjoining all
primes $pL$ over $p,$ for all $L,$ along with $\pr^{-1}(pL).$   
  
\bigskip

{\bf Remark.} Suppose $\cl(X)\cap \cl(Y)=\emptyset,$ for $X,Y\subset   ({^*\mathcal{AF}}_K, {^*\mathcal{PF}}_K).$
Since the language does not contain predicates and functions  linking $\cl(X)$ and $\cl(Y),$ for any disjoint embeddings $e_X$ and $e_Y,$   $({^*\mathcal{AF}}_K, {^*\mathcal{PF}}_K)\to ({\mathcal{AF}'}_K, {\mathcal{PF}'}_K),$  
of   $\cl(X)$ and $\cl(Y),$ the union  $e_X\cup e_Y$ is an embedding.

\subsection{Lemma}\label{stableclO} {\em   In $({^*A}_K, {^*\Spec}_K),$ for any point $a,$ the substructure on the set $\cl(a)$ is $\omega$-stable of finite Morley rank on each sort.

Moreover,  any formula $\varphi(v_1,\ldots,v_n)$ relativised to   $\cl(a)$ is equivalent to an $\exists$-formula.
}

{\bf Proof.} As described in the Example in \ref{cl}, the infinite definable subsets consist of $O_K$-monoid-modules  (orbits) $\mathrm{Fb}_{pK}$ with monoid homomorphism $\pi_{K,L}:   \mathrm{Fb}_{pK}\to \mathrm{Fb}_{pL}$ between those. This is $\omega$-stable of finite Morley rank.

\begin{comment}
To prove the second statement note that any formula $\varphi(v_1,\ldots,v_n)$ involves variables restricted to finitely many sorts $A_K$ and $\Spec_K.$ If it is also relativised to $\cl(a)=\cl(p),$ ($p\in \Spec_\Q$) then the formula is relativised to a finite structure. This finite structure  is completely determined by the sorts involved and $p.$
 \end{comment}
 
The second statement is equivalent to the fact that the theory of $\cl(a)$ is model-complete. To prove the latter consider the structures in two embedded models, $\cl(a)\subseteq \cl(a)'.$ It is clear from the description in \ref{cl} that $\cl(a)$ and $\cl(a)'$ have the same primes $pK$ in all sorts. Hence the embedding reduces to the family of embeddings  $\mathrm{Fb}_{pK}\subseteq  \mathrm{Fb}'_{pK}.$  The Tarski-Vaught test verifies that this is an elementary embedding,  $\cl(a)\preceq \cl(a)'.$
\qed

\subsection{Lemma}\label{stableclF} {\em   In  $({^*\mathcal{AF}}_K, {^*\mathcal{PF}}_K),$ for any point $a,$ the substructure on the set $\cl(a)$ is supersimple of finite rank on each sort.}

{\bf Proof.} As described in the Example in \ref{cl}, the infinite definable subsets consist of fields $\F_{pK}$
and its definable subfields $\F_{pL}$. Each of these is a finite extension of $\F_{p{\Q}},$ the minimal field in $\cl(a),$
and so interpretable in $\F_{p{\Q}}.$
The latter is supersimple of finite rank. \qed

  \subsection{Lemma}\label{3.7} {\em Assume the continuum hypothesis,  let  $\mathcal{D}$ of \ref{UfD} be a non-principal good ultrafilter on a countable set of indices and let  $({\mathcal{AF}'}_K, {\mathcal{PF}'}_K)$ be any saturated  space of pseudofinite fields  of cardinality $\aleph_1.$ Then for any algebraically closed $M\subset ({^*\mathcal{AF}}_K, {^*\mathcal{PF}}_K)$ and algebraically closed $M'\subset ({\mathcal{AF}'}_K, {\mathcal{PF}'}_K)$ 
   such that $\Spec_L\cap M$ is countable for any $L,$ and $M\cong_\mathbf{e}M'$ % by an isomorphism $\mathbf{e},$
   there is an extension of the  isomorphism $\mathbf{e}$ to an isomorphism  
   $$({^*\mathcal{AF}}_K, {^*\mathcal{PF}}_K)\cong ({\mathcal{AF}'}_K, {\mathcal{PF}'}_K).$$}
 
 {\bf Proof.} Note that by our assumptions $({^*\mathcal{AF}}_K, {^*\mathcal{PF}}_K)$ is also a saturated  space of pseudofinite fields  of cardinality $\aleph_1.$  So we need to construct an isomorphism between any two of those.
 
 It is enough to show that for every $a\in  {^*\mathcal{AF}}_K\setminus M$ there is an  $a'\in  {\mathcal{AF}'}_K\setminus M'$ such that $\mathbf{e}$ can be extended to  an isomorphism $M\cup \cl(a)\to M'\cup \cl(a').$ 
 
 By assumptions $a\in \F_q$ for some pseudofinite field $\F_q,$ for $q=q(a)\in \Spec_K,$ some number field $K.$
 The type $\mathrm{tp}(a,\F_q)$
 of $\F_q$ and $x$ in the language    $L_\mathrm{rings}(K)$ determines, by axioms (\ref{axiomFqq})--(\ref{Nm}) and \ref{cl} the type    $\mathrm{tp}(a,\cl(a)).$ 
 
 On the other hand, since $\F_q$ and all the other fields in $\cl(a)$ are pseudofinite, each formula in $\mathrm{tp}(a,\F_q)$ is realised by some pair with a finite field $\F_{p^n}$ in place of $\F_q.$ It implies that the same collection of formulas is a type in $({\mathcal{AF}'}_K, {\mathcal{PF}'}_K),$ and so has a realisation $a', \F_{q'}.$ Hence, also
 $$\mathrm{tp}(a,\cl(a))=\mathrm{tp}(a',\cl(a')).$$
Since both $\cl(x)$ and $\cl(x')$ can be seen as sort-by-sort definable substructures of saturated structures, we have an isomorphism   $\cl(a)\to \cl(a'),$ $a\mapsto a'.$  By the Remark in \ref{cl} we obtain 
an isomorphism $M\cup \cl(a)\to M'\cup \cl(a')$ extending $\mathbf{e}.$ \qed

\subsection{Corollary}\label{3.8} {\em Any two spaces of pseudofinite fields are elementarily equivalent. The type of  the ``flag'' $\langle a, q\rangle,$ $a\in \F_q,$ over an algebraically closed $M,$ $q\notin M,$ is determined by the elementary theory of the pair in the language  $L_\mathrm{rings}(K).$ The complete invariant of the theory of $\F_q$ is the size of $\F_q,$ if $\F_q$ is finite,  or, otherwise, the subfield $\F_q\cap \tilde{\Q},$ with names from $O_K.$ }

For the last part of the statement see \cite{Zoe}.

\subsection{Lemma.}\label{new} {\em In the structure  $(^*A_K,^*\Spec_K)$  let $q$ be an  infinite prime and $\F_{q}$  a pseudo-finite field containing $K$ as the named subfield. Let $M \supset K$, then the following are equivalent.

(i) $\F_{q}$ contains some Galois conjugate $M'\supseteq K$ of $M$ as a named field;

(ii) for some Galois conjugate $M'\supseteq K$ of $M$ and   prime $ q' $ of $\mathrm{O}_{M'}$ lying above $q$, $\F_{q^{'}}\cong\F_{q}.$}

%(iii) $q\in \bigcup_{K'\cong_L K}\Pi_{K',L}$ (all Galois conjugates $L'$ of $L$).}

\medskip

\begin{proof} Write ${^*\mathrm{O}}_{M'}$ for the ultrapower of $\mathrm{O}_{M'}.$ This will by definition contain $q'$ as an ultraproduct of prime ideals.

$(ii) \Rightarrow (i)$, since $\mathrm{res}_{q'}(^*\mathrm{O}_{M'})={^*\mathrm{O}}_{M'}/q' = \F_{q'}\cong \F_{q}$ and $\mathrm{res}_{q'}$ is injective on $M'.$

$(i) \Rightarrow (ii)$ Suppose $\mathrm{O}_{M'}$ has a basis $b_1,\cdots b_n$ as an $\mathrm{O}_{K}$-module. 
Then $b_1,\cdots,b_n$ is also a basis for $^*\mathrm{O}_{M'}$ over $^*\mathrm{O}_{K}.$ 

Now if $(i)$ holds,  
let $q'$ be a prime of ${^*\mathrm{O}}_{M}$ lying above $q$, the residue fields $^*\mathrm{O}_{M}/q'=\F_{q'}$ would be a finite extension of $^*\mathrm{O}_{K}/q=\F_{q}$ generated by  the named elements $\mathrm{res}_{q'}(b_1),\cdots,\mathrm{res}_{q'}( b_n)$ which are in $\F_{q}$ by assumption. Then $\F_{q'}= {^*\mathrm{O}_{M}}/q'={^*\mathrm{O}_{K}}/q=\F_{q}.$
%, (in case of finite fields if $\F_{q}$ contains a conjugate of $\F_{q'}$, implies that $\F_{q}=\F_{q'}$).

%$(ii) \Leftrightarrow (iii)$ by definition \ref{3.14+}.
\end{proof}

\subsection{Corollary}\label{3.9}{\em Let $N$ be an algebraically closed substructure of the multisorted representations structure $({^*A}_K,{^*\Spec}_K),$ and $q,q'\in \Spec_K.$ Then 
\begin{equation}\label{3.9a}
\mathrm{tp}(q/N)=\mathrm{tp}(q'/N) \Leftrightarrow \F_q\equiv_{L_\mathrm{rings}(K)}\F_{q'}
\end{equation}

and, for $q$ infinite,
 \begin{equation}\label{3.9b}
\mathrm{tp}(q/N)=\mathrm{tp}(q'/N) \Leftrightarrow \mbox{ for all } M\supset K:\ q\in \Pi_{M,K} \leftrightarrow  q'\in \Pi_{M,K}\end{equation}
 }

%Given $a_1,\ldots, a_n\in O_q$ the type $\mathrm{tp}(a_1,\ldots, a_n, q/N)$ is determined by $\mathrm{tp}(q/N)$ and the collection of formulas of the form $\gamma_{ij}\cdot a_i=a_j,$ for 
%$\gamma_{ij}\in O_K$ for some $i,j\in \{1,\ldots, n\}.$
\begin{proof}

(\ref{3.9a}) follows from \ref{UfD},  \ref{3.7} and \ref{3.8}.

In order to prove (\ref{3.9b}) invoke 
Kiefe's criterium (see \cite{Zoe}, 4.7):   $ \F_{pK}\equiv \F_{p'{K}}$ in the language $L_\mathrm{rings}(K)$ if and only if for every irreducible polynomial $f(x)$ over $K,$ 
$$\F_{pK}\vDash \exists x\, f(x)=0 \Leftrightarrow \F_{p'K}\vDash \exists x\, f(x)=0,$$
where coefficients of $f(x)$ are interpreted in $\F_{pK}$ and $\F_{p'K}$ as names of elements in the pseudofinite fields.  

Note that the assumption that $pK\in \Spec_K$ is infinite (i.e. non-standard) implies that $\mathrm{char}\,\F_{pK}=0$ and the naming homomorphism $K\to \F_{pK}$ is an embedding. It follows that $f(x)$ remains irreducible as a polynomial over $K_{pK}=K,$ the image of $K$ in $\F_{pK}.$

Let $M\supset K$ be a field generated by a root of $f(x)$ over $K.$ Now 
$$\F_{pK}\vDash \exists x\, f(x)=0 \Leftrightarrow \bigvee_{M'\cong_K M}\F_{pK}\supseteq M'$$
where on the right we consider images of conjugates of $M$ in $\F_{pK}.$
But this condition by \ref{new} and \ref{3.14++}  is equivalent to $pK\in \Pi_{M,K}.$
\end{proof}

\subsection{Remark}
\label{rem4.12} (i) {\em Suppose $M/ K$ is Galois. Then for all but finitely many $q\in \Spec_K,$
$$q\in \Psi_{M,K}\Leftrightarrow q\in \Pi_{M,K}.$$}
Indeed, let  
 as above $M$ be generated by a root of an irreducible $f(x)$ of order $n$ over $K.$ Then $M$ contains all the roots of $f$ and so for an infinite prime $q$
  $$\F_{q}\vDash \exists x\, f(x)=0 \Leftrightarrow \F_{q}\supseteq M \Leftrightarrow 
\F_{q}\vDash \exists x_1,\ldots, x_n\,  \left( \bigwedge_{i=1}^n f(x_i)=0\ \& \ \bigwedge_{i\neq j} x_i\neq x_j   \right).
  $$
  The condition on the right means that $f(x)$ splits completely over $q,$ that is $q\in  \Psi_{M,K}.$ The condition on the left means that $q\in  \Pi_{M,K}.$ This holds for all nonstandard primes $q, $ hence does hold for all but finitely many standard primes. 

(ii)  {\em The condition that $M/K$ is Galois is not essential.   Let $\hat{M}$ be the minimal Galois extension of $K$ containing $M$ Then
 for all but finitely many $q\in \Spec_K,$
$$q\in \Psi_{M,K}\Leftrightarrow q\in \Psi_{\hat{M},K}\Leftrightarrow q\in \Pi_{\hat{M},K}.$$}
Indeed, under the assumtion we have an infinite $q\in \Spec_K$ such that $q\in  \Psi_{M,K}.$ 
But then by definition of $ \Psi_{M,K}$ we will have that the minimal polynomial $f(x)$ for $M$ splits completely (see \ref{3.14+}). Then $\F_q$ contains all roots of $f(x),$ so all  conjugates of $M$ so $\hat{M}.$ Thus $q\in  \Psi_{\hat{M},K}.$ Conversely, if  $q\in  \Psi_{\hat{M},K}$ then
$\F_q$ contains $\hat{M},$ so contains all conjugates of $M$ and so $f(x)$ splits completely in $\F_q,$ $q\in   \Psi_{M,K}.$ 
 %which generated $M$ over $K,$ and so $M$ is invariant  

 \subsection{Corollary (to \ref{3.9})} {\em The theory of the multisorted representations structure is superstable. The $U$-rank of a non-principal 1-type of sort $\Spec_K$ is $1,$ the $U$-rank of a 1-type  of $A_K$ containing the formula $\pr(v)=p$ over $p\in \Spec_K$ is $0$ or $1$ depending on whether $p$ is finite or infinite prime. If the type contains no such formula, it is generic and its $U$-rank is 2.}
 
 \medskip
 
 Indeed, by (\ref{3.9b}) for any $N$ the set of $1$-types over $N$ of sorts  $\Spec_K$ is at most of cardinality continuum. Moreover, each $0$-definable non-principal type has a unique non-algebraic extension over  $N.$
   
  A type of sort $A_K$ by \ref{3.7}  is determined by a type of sort $\Spec_K$ and the type relative to the $\omega$-stable substructure $\cl(a).$ If a 1-type contains a formula $\pr(v)=p$ for a parameter $p,$ then it is a type of an element of the fibre $\mathrm{Fb}_p,$ which is either finite, if $p$ is finite, or strongly minimal, so of $U$-rank 1.  The unique complete 1-type which negates all formulas $\pr(v)=p$ is of $U$-rank 2.\qed

 \subsection{Theorem}\label{ThmQE}
{\em The theory of the multisorted representations structure has QE in the language extended by boolean combination of the unary predicates:
\begin{itemize}
\item $v\in \Spec_K,$
\item $v=q,$ for $q\in \Spec_K,$ finite,
\item  $v\in \Pi_{K,L}$
\end{itemize}
  and 
\begin{itemize}
\item  existential formulas  $\varphi(v_1, w_1,\ldots,v_n,w_n,p)$ where the $v_i$ are of sorts $A_{L_i}$ and $w_i$ of $\Spec_{L_i}$ and $v_1,w_1\ldots, v_n,w_n$ are
 relativised to   $\cl(p),$  $p\in \Spec_\Q,$ $L_1,\ldots, L_n\in \mathcal{R}.$
\end{itemize}  
 }

\smallskip

{\bf Proof.} By \ref{stableclO} and  \ref{3.9}(\ref{3.9b})  a complete $n$-type in the theory  is determined by formulas listed in the formulation of the theorem. By the compactness theorem any formula is equivalent to a boolean combination of those.\qed

\section{The topology and compactification}

We extend the conventional topology of \ref{top} on  $\Spec_{L}$ 
by declaring {\bf closed in $\Spec_{L}$}
subsets of the form $\Pi_{K,L},$ singletons  and their finite unions, for all  extension $K$ of $L.$ 

 We accordingly extend the topology on $A_{L}$ by declaring {\bf closed  in  $A$} 
 \begin{itemize}
 \item the graph of $\pr$ restricted to $A_L\times \Spec_L,$  $L\in \mathcal{R};$
 \item  the graph of $\pi_{K,L},$ for   $K, L\in \mathcal{R};$
 \item   subsets of  $\Spec_L$ of the form $\Pi_{K,L}$, 
 $L\in \mathcal{R};$
 \item subsets of $A_{L_1}\times\Spec_{L_1} \ldots \times A_{L_n}\times \Spec_{L_n}$  defined by positive existential formulas  $\varphi(v_1, w_1,\ldots,v_n,w_n,p)$ relativised to   $\cl(p),$  $p\in \Spec_\Q,$ $L_1,\ldots, L_n\in \mathcal{R}.$
 \end{itemize}
along with cartesian products, finite intersections and unions of those.

\subsection{Lemma} {\em  In the multisorted structure
% with sorts associated to Galois extension $\Q[\alpha]$ of $\Q,$ 
the maps $\pr$ and $\pi_{K,L}$ are continuous. }

{\bf Proof.} The continuity of $\pr$ is just by definition.

In order to prove the continuity of $\pi_{K,L}$ we need to prove that $\pi^{-1}_{K,L}(S)$ is closed for every closed $S\subset A_L.$ Since this is obvious for finite $S,$ we need to
 consider only the closed subsets $S$ of the form $\Pi_{M,L}.$

{\bf Claim.} $$\pi^{-1}_{K,L}(\Pi_{M,L})=\Pi_{KM,K}$$
where $KM$ is the composite of fields.

Proof.
Let $pL\in\Pi_{M,L},$ $M=L[\alpha].$ 
Then we have $KM=K[\alpha]$ and for some $a\in L,$
\begin{equation} \label{*}
\alpha-a\in pL.
\end{equation} 
 
Then  
\begin{equation} \label{**}
\alpha-a\in pK\mbox{ for every prime } O_K\supset pK\supseteq pL
\end{equation} 
which implies $pK\in \Pi_{KM,K}.$ 
This proves   
 $$pL\in\Pi_{M,L}\Rightarrow pK\in \Pi_{KM,K}\mbox{ for every prime } O_K\supset pK\supseteq pL.$$
 
Now we prove the converse. Suppose  
$$ pK\in \Pi_{KM,K}\mbox{ for every prime }  pK\supseteq pL.$$
This is equivalent to say that $\F_{p{KM}}=\F_{pK}$ and so to
(\ref{**}).
But $$\bigcap_{pK\supseteq pL} pK= pL$$ and hence we proved 
(\ref{*}) and  $pL\in\Pi_{M,L}.$ This finishes the proof of the claim and of the lemma. \qed

\subsection{The Chebotarev density theorem and its corollaries} \label{chebotarev}
%\footnote{Need to do it for any, non-Galois, extension, as well as Thm \ref{3.18}.}
 Recall that the {\bf density }of a subset  $S\subset \Spec_{L}$ is defined as 
 $$\mathrm{dn}\,S=\lim_{n\to \infty} \frac{\#\{ \mathfrak{p}\in S: |\mathfrak{p}|\le n\} }{\#\{ \mathfrak{p}\in \Spec_{L}: |\mathfrak{p}|\le n\}}.$$
 
The  Chebotarev density theorem (\cite{MP}, 4.4.3) states that a subset  $S\subset \Spec_{L}$  defined by the pattern of splitting of $f_K,$ the   minimal  polynomial for $K$ over $L,$ modulo $\mathfrak{p}$ has a density. Moreover, it determines the density in terms of the structure of $G=\mathrm{Gal}(\hat{K}:L),$ where $\hat{K}$ is the minimal Galois extension of $L$ containing $K.$   More precisely, let  $\sigma\in G$ and
$C=\sigma^G$ a conjugacy class.

Then Chebotarev Density Theorem states that
$$\mathrm{dn} \{ p\in \Spec_L:\ 
\sigma_p\in C \}=\frac{|C|}{|G|},$$
where $\sigma_p$ are the Frobenius elements defined up to conjugacy.

\subsection{Proposition}\label{3.17} {\em 
Let $L\subset K\subset M$ and $M/L$ minimal Galois extension containg $K.$  

 Then $\mathrm{dn}(\Pi_{M,L}) $ is well-defined and 
$$\mathrm{dn}(\Pi_{K,L}) \ge \frac{1}{|M:L|}.$$}
%{\mbox{no of Galois conjugate of } L \mbox{ in } K}$$

\begin{proof}

%There are only finitely many $p\notin \Delta_{M/K}$ so we may ignore the first condition on $p.$

Let $p=pL,$ $q=pM$.
Consider the Frobenius element $\mathrm{Frob}_p\in G=\mathrm{Gal}(M:L),$  that is such that its action on $\F_q=\mathrm{res}_q(M)$ 
 fixes exactly $\F_p$ inside $\F_q.$  Now 
$\mathrm{res}_q(K)=\F_r,$   the image of $K$ in $\F_q,$ is fixed by the action of $\mathrm{Frob}_p^m,$ where $m=|\F_r:\F_p.|$   

Set $C\subseteq G$ to be the conjugacy class which contains $\mathrm{Frob}_p^m.$  
Clearly $|C|\ge 1$ and $p$ is such that   
the Chebotarev formula above now counts the density of those $p\in \Spec_L$ for which there
is (every) $q\in \Spec_M,$ $q\supset p,$     $\mathrm{Frob}_p^m=\mathrm{Frob}_p.$ That is 
$\F_r=\F_p,$ equivalently $\F_{pL}\cong \F_{pK}.$ These are exactly the $pK$ in  $\Pi_{K,L}.$
%$\mathrm{O}_L/pL\cong \mathrm{O}_K/pK.$
\end{proof}

 An immediate corollary to \ref{3.17} is that    $\Pi_{K,L}$ and $\Psi_{K,L}$ are infinite. Also
 $\Psi_{K,L}$ has well-defined density,
 $$\mathrm{dn}(\Psi_{K,L}) = \frac{1}{|\hat{K}:L|},$$
 where $\hat{K}$ is the minimal Galois extension of $L$ containing $K.$ See \ref{rem4.12}.

\subsection{Lemma}\label{intersect0} {\em 
For any two  Galois extensions $K_1$ and $K_2$ of $L$ and any infinite prime $q\in \Spec_L,$
%the splitting fields for 
%$f_{K_1}$ and $f_{K_2}$ modulo $\mathfrak{p}$ is the composite of the splitting fields for each of those, and is equal to the splitting field of $f_K,$ for $K=\Q[\alpha_1,\alpha_2],$ we get that for Galois extensions
$$q\in \Psi_{K_1,L}\cap \Psi_{K_2,L}\Leftrightarrow q\in \Psi_{K,L}$$
for $K=K_1K_2,$ the composite of the two fields
 and so any intersection of the $\Psi_{K_i,L}$ has a well-defined positive density.} 

{\bf Proof.}  Indeed, $K_1$ and $K_2$ are subfields of $\F_q$ if and only if $K$ is. $\Box$

\subsection{Proposition} \label{intersect}
An intersection of finite family of special subsets  $\Pi_{K,L}$ (arbitrary $K:L$)  of $\Spec_L$ is nonempty:
 $$\bigcap_{i\in I}\Pi_{K_i,L}\neq \emptyset.$$
 
 {\bf Proof.} First note that by definition
 $$\bigcap_{i\in I}\Pi_{K_i,L}\supseteq \bigcap_{i\in I}\Psi_{K_i,L}.$$

 By   \ref{intersect0} and \ref{rem4.12},   for infinite $q\in \Spec_L$ 
  $$q\in \bigcap_{i\in I}\Psi_{\hat{K}_i,L}\Leftrightarrow q\in \Psi_{K,L}\Leftrightarrow q\in \Pi_{K,L},$$
 for $K$ the composite of all the  $\hat{K}_i.$ Hence  $\bigcap_{i\in I}\Psi_{\hat{K}_i,L}$ is infinite and non-empty.
$\Box$

\medskip

We note also the following.
\subsection{Theorem}\label{3.18}
 {\em For any  extensions $K_i$ of $L,$ $i=1,\ldots, k,$ 
any boolean combination of subsets $\Psi_{K_i,L}$ of $\Spec_{L}$ 
has a well-defined density.
%\begin{equation}\label{intersect}
%\mathrm{dn}\,(\bigcap_{i=1}^k\Pi_{\alpha_i,\beta})> 0 
%\end{equation}
%In particular, the intersection is infinite.
}

{\bf Proof.} First we note that by \ref{rem4.12}(ii) we may assume that $K_i/L$ are Galois.

Secondly, note that 
 it is enough to prove the theorem for the intersection  
$$\bigcap_{i\in \{ i_1,\ldots i_m\}}\Psi_{K_i,L}\ \cap \
\bigcap_{j\in \{i_{m+1},\ldots i_k\}} \neg  \Psi_{K_j,L}$$
for any permutation $\{ i_1,\ldots,i_k\}$ of $\{ 1,\ldots,k\}.$ 

Also, as shown in \ref{intersect0} $\bigcap_{i\in \{ i_1,\ldots i_m\}}\Psi_{K_i,L}$ can be replaced by some
$\Psi_{K,L}.$

Now we prove,

\smallskip

{\bf Claim}. The density of the set of the form $$\bigcup_{l\in \{ 1,\ldots m\}}\Psi_{K_l,L}$$ is well-defined.

We prove this by induction on $m.$ Assuming it is true for $m$
note that by definition
$$\dn \left(\bigcup_{l\in \{ 1,\ldots m+1\}}\Psi_{K_l,L}\right)=
\dn \left(\bigcup_{l\in \{ 1,\ldots m\}}\Psi_{K_l,L}\right) +\dn \, \Psi_{K_{m+1},L} - \dn\, \left((\bigcup_{l\in \{ 1,\ldots m\}}\Psi_{K_l,L})\cap \Psi_{K_{m+1},L}\right)$$
provided the three summands on the right are well-defined. The first two are well-defined by the induction hypothesis. The same is true for the last one since
$$\left(\bigcup_{l\in \{ 1,\ldots m\}}\Psi_{K_l,L}\right)\cap \Psi_{K_{m+1},L}=\bigcup_{l\in \{ 1,\ldots m\}}\left(\Psi_{K_l,L}\cap \Psi_{K_{m+1},L}\right).$$
This proves the claim.
\medskip

Now it remains to see that 
$$\dn \left(\Psi_{K,L} \cap 
\bigcap_{j\in \{1,\ldots m\}} \neg  \Psi_{K_j,L}\right)=
\dn \Psi_{K,L}\ - \ \dn
\bigcup_{j\in \{ 1,\ldots m\}}  \Psi_{K,L}\cap \Psi_{K_j,L}.
$$
The last term is well-defined by the claim. This finishes the proof of the theorem.\qed

\medskip

{\bf Remark.}  The density is a measure on the boolean algebra generated by all the $\Psi_{K,L}$ for a given $L.$

\medskip

{\bf Question.} Is  density well-defined on the boolean algebra generated by all the $\Pi_{K,L},$ for a given $L?$

\subsection{Compact models}
Note that the sorts  $\Spec_{L}$ and $A_{L}$ are not compact since the intersection 
$\bigcap_{i\in I}\Pi_{K_i,L},$
for an infinite $I$ and distinct $K_i,$ $i\in I,$ is empty (no prime ideal in $O_{L}$ belongs to such an intersection) but the intersection of any finite subfamily of the sets is non-empty  by \ref{intersect}.

However, there are plenty of models which are compact in the
special topology; all we need to do is to realise the maximal positive types  \begin{equation}\label{limitp0}
\bigcap_{K\supset L}\Pi_{K,L}
\end{equation}
in each sort $\Spec_L$ (these are types by \ref{intersect}).

These contain positive types of the form \begin{equation}\label{limitp}
\bigcap_{K\supset L}\Psi_{K,L}
\end{equation}
 Any realisation of such type will be called a {\bf splitting infinite prime} (or just a splitting prime) in $\Spec_L.$

We call a multisorted structure $\mathfrak{M}$
of representations {\bf a minimal compact model} if 
\begin{itemize}
\item[(i)]
 $\mathfrak{M}$ is an elementary extension of the standard  multisorted structure 
of representations;
\item[(ii)]
any sort $\Spec_L(\mathfrak{M})$ in the structure is compact;
\item[(iii)] for any other  $\mathfrak{M}'$ satisfying (i) and (ii), for each sort $A_L$ there is an embedding 
 
 $$A_L(\mathfrak{M})\subset A_L(\mathfrak{M}').$$

\end{itemize}
\subsection{Theorem}\label{compact}
 {\em A minimal compact model $\mathfrak{M}$ %of the multisorted theory of representations
exists and is unique.

(a)  For each Galois extension $L$ of $\Q$ a minimal $\mathfrak{M}$ contains exactly $\deg L$ infinite primes in $\Spec_L.$ 

(b) All the infinite primes in $\Spec_K(\mathfrak{M})$ are splitting primes and all the infinite primes containing 
$pL\in \Spec_L(\mathfrak{M})$  for $K:L$ Galois are conjugated by the definable action of the Galois group $\mathrm{Gal}(K:L)$ defined by formula  (\ref{qq}). 

(c) The residue field  $\F_{p{L}}$ over such an infinite prime $pL$ is characterised as a pseudo-finite field containing $\Q^{alg}$ with a naming homomorphism 
  $$\mathrm{res}_{p{L}}: \mathrm{O}_{L}\to\Q^{alg}\subset \F_{p{L}}$$
 The first order theory of such  $\F_{p{L}}$ is determined uniquely by the latter inclusion.

(d) The fibres over infinite $pL$ have the form
 $$\mathrm{Fb}_{pL}(\mathfrak{M})=\Q^{alg}\cdot a_{pL}$$
 for a non-zero element $a_{pL}$ of the fibre.  
}    
\smallskip 
 
{\bf Proof.} We use the analysis in \ref{psi} for finite primes in $\Pi_{K,L}$ and extend the conclusions to infinite primes satisfying the same first-order condition. The Tarski-Vaught test along with the statement of theorem \ref{ThmQE} allows to conclude that the structure with properties (a)-(c) is an elementary submodel of a saturated model $({^*A}_K, {^*\Spec_K}),$ $K\in \mathcal{R},$ described in section \ref{QE}.  

It remains to establish properties (a)-(d).

\medskip

(a)-(b). In our  case the splitting occurs  that is $\F_{p{K}}=\F_{p{L}}$ and $\mathrm{Norm}_{K,L}$ is the identity map.  The  $|K:L|$  distinct  prime ideals $p^{(i)}{K} \in \Spec_{K},$ ($i=1,\ldots, |K:L|$) lying over $p{L}$  give rise to
$|K:L|$  naming  homomorphisms 
 $$\mathrm{res}_{p^{(i)}{K}}: \mathrm{O}_{K}\to \F_{p{L}},\ \ \ \mathrm{ker}(\mathrm{res}_{p^{(i)}{K}})={p^{(i)}{K}}$$
 which are pairwise non-isomorphic over $\mathrm{res}_{p{L}}: \mathrm{O}_{L}\to \F_{p{L}}.$
 
Hence we have $|K:L|$ different ways of naming elements of $\F_{p{L}}$ when we pass from level
$L$ to level $K.$ This proves that over infinite splitting $p\in \Spec_\Q$ there are exactly $\deg L$ splitting $pL$ in $\Spec_L.$

The distinct prime ideals $p^{(i)}{K}$ extending $pL$ in $O_K$ are conjugated under the action of 
 $\mathrm{Gal}(K: L),$ which is first-order definable by (\ref{qq}), so extends to the infinite splitting primes,  and this proves the second statement.

(c) follows from (\ref{3.9a}) and
the characterisation of the elementary types of pseudo-finite fields $\F$ of characteristic $0$ by their subfield of algebraic numbers (see \cite{Zoe}).

(d) The only condition on the infinite fibres in $A_K$ is that the action by $K$ is defined and is free. This is satisfied if we set the fibre of the form (d).    $\Box$
   
\subsection{The minimal complete model} We call  $\mathcal{N}$
a  minimal complete model  if:  
\begin{itemize}
\item[(i)]
 $\mathfrak{N}$ is an elementary extension of the standard  multisorted structure 
of representations;
\item[(ii)]
any sort $\Spec_L(\mathfrak{N})$ in the structure realises all the 1-types over $0$;
\item[(iii)] for any other  $\mathfrak{N}'$ satisfying (i) and (ii), for each sort $A_L$ there is an embedding 
 
 $$A_L(\mathfrak{N})\subset A_L(\mathfrak{N}').$$

\end{itemize} 

\subsection{Theorem}\label{complMod} {\em A minimal complete model exists and is unique. 

(a)  For each Galois extension $L$ of $\Q$ a minimal $\mathfrak{M}$ contains at most $\deg L$ infinite primes in $\Spec_L.$ The primes $pK$ over $pL,$ for $K/L$ Galois, are conjugated under the action of Galois group $\mathrm{Gal}(K:L).$

(b) An infinite prime $pL\in \Spec_L(\mathfrak{N})$ is uniquely characterised by the field $$E_{pL}=\Q^{alg}\cap \F_{p{L}}$$ and the naming homomorphism 
  $$\mathrm{res}_{p{L}}: O_L\to E_{pL}.$$

(c) The fibres over infinite $pL$ have the form
 $$\mathrm{Fb}_{pL}(\mathfrak{N})=E_{pL}\cdot a_{pL}$$
 for a non-zero element $a_{pL}$ of the fibre.   
}

\smallskip

{\bf Proof.} Same arguments as in \ref{compact}.$\Box$

\subsection{Formal geometry} We leave out the problem of identifying a model which can be seen as a generalised Zariski geometry.
It is clear that neither of the models we discussed above is an analytic (or Noetherian) Zariski geometry in the sense of \cite{ZBook}.

\section{The  universal cover}
% of $\mathfrak{M}_{\mathrm{Rep}}$}

\subsection{The coarse cover.} Let $\mbox{O}_\mathcal{R}=\bigcup\{ \mbox{O}_{K}: \mbox{O}_{K}\in \mathcal{R}\}.$
Most of the time we will deal with the case of $ \mathcal{R}$ containing all integral extensions of $O_M,$ some number field $M.$

The systems of morphisms $\pi_{K,L}$ and  $\pi^\Spec_{K,L},$ $K\supset L,$ determine a projective system. 

Define $A_\mathcal{R}$ and 
$\Spec_{\mathcal{R}}$ 
as the projective limits of  $A_K$ and 
$\Spec_{K}$ respectively,   $\mbox{O}_{K}\in \mathcal{R}.$

By definition we will have 
maps 
$$\pi_{\mathcal{R},K} :A_{\mathcal{R}} \to A_{K},\ \mbox{ and } \ \pi^\Spec_{\mathcal{R},K} :\Spec_{\mathcal{R}} \to \Spec_{K}$$
satisfying  %$$\pi_{\mathcal{R},\alpha} \circ \pi_{\alpha,\beta}= \pi_{\mathcal{R}\beta}\mbox{ and }
$$\pi_{\mathcal{R},K} \circ \pi_{K,L}= \pi_{\mathcal{R},L},\ \mbox{ and } \ \pi^\Spec_{\mathcal{R},K} \circ \pi^\Spec_{K,L}= \pi^\Spec_{\mathcal{R},L}.$$ 

Also the following (non-first-order) condition is satisfied:
\begin{equation} \forall x_1,x_2\in A_\mathcal{R}\ x_1=x_2 \leftrightarrow \bigwedge_{K} \pi_{\mathcal{R},K}(x_1)=\pi_{\mathcal{R},K}(x_2)
\end{equation}

The definition assigns to each point $\tilde{p}\in  \Spec_{\mathcal{R}}$ an 
 orbit $\mathrm{Fb}_{\tilde{p}}$
which, 
 according to \ref{pi}, can be represented  as
 $\mathrm{Fb}_{\tilde{p}}=\F_{\tilde{p}}.a_{\tilde{p}}$ for some $a_{\tilde{p}}\in \mathrm{Fb}_{\tilde{p}},$ where by definition 
 $$\F_{\tilde{p}}=\bigcup\{ \F_{p{K}}:  \pi_{\mathcal{R},K}(\mathrm{Fb}_{\tilde{p}})=\F_{p{K}}.a_ {\pi_{\mathcal{R},K}(\tilde{p})}\},$$
the union of the tower of named fields below $\tilde{p}.$ The naming means that $\F_{\tilde{p}}$ is given along with 
a naming homomorphism $$\mathrm{res}_{\tilde{p}}:{\mbox{O}}_\mathcal{R}\to  \F_{\tilde{p}}.$$

Clearly,  $$A_\mathcal{R}:=\bigcup\{ \mathrm{Fb}_{\tilde{p}}:  \tilde{p}\in  \Spec_\mathcal{R} \}.$$

The topology on $A_\mathcal{R}$ and $\Spec_\mathcal{R}$ is defined as the projective limit of topologies on the $A_K$ and $\Spec_K,$ that is a subset $S\subseteq \Spec_\mathcal{R}$ is defined to be closed 
if $\pi_{\mathcal{R},K}(S)$ is closed for all large enough $K.$

In particular, the fibre
$(\pi^\Spec)^{-1}_{\mathcal{R},K}(q)\subset \Spec_\mathcal{R}$ is closed for any $K$ and $q\in \Spec_K.$ 

\subsection{The fibre of $A_\mathcal{R}$ over a finite prime}
In case when  $\tilde{p}$ lies over a standard prime $p$, that is 
$\pi^\Spec_{\mathcal{R},1}(\tilde{p})=p\in \Z,$ and $\bigcup \mathcal{R}=\Q^{alg}$  we have  $\F_{\tilde{p}}=\F_p^{alg}.$ However, for the same $p$ we will have as many points $\tilde{p}$ over $p$ as there are  naming homomorphisms ${\mbox{O}}_ \mathcal{R}\to  \F_{\tilde{p}}.$ 
Note that the naming homomorphisms in this case are just residue maps of $p$-adic valuation on $\Q^{alg}.$
In other words,  setwise
 $$(\pi^\Spec_{\mathcal{R},1})^{-1}(p)=\{ p\mbox{-adic valuations on }\Q^{alg} \}\cong \{ \mathrm{res}_p: {\mbox{O}}_\mathcal{R}\to  \F_p^{alg}\}.$$

 Let $\mathcal{G}_\Q:=\mathrm{Gal}(\tilde{\Q}:\Q).$ Then, it follows from the general theory of valuations that, for any two homomorphisms $\mathrm{res}_p$ and $\mathrm{res}'_p:
 {\mbox{O}}_\mathcal{R}\to  \F_{\tilde{p}},$ there is a $\sigma\in \mathcal{G}_\Q$ such that $\mathrm{res}'_p=\mathrm{res}_p\circ \sigma.$ In other words $\mathcal{G}_\Q$ acts transitively on the fibre
and  $$(\pi^\Spec_{\mathcal{R},1})^{-1}(p)=\tilde{p}\cdot \mathcal{G}_\Q.$$ 
%for $\tilde{p}\in (\pi^\Spec_{\mathcal{R},1})^{-1}(p).$

Note that $$\mathrm{ker}(\mathrm{res}_p\circ \sigma)=\sigma^{-1}\{ \mathrm{ker}(\mathrm{res}_p)\}$$
and so the subgroup of  $\mathcal{G}_\Q$ fixing $\tilde{p}$  is equal to the stabiliser $\mathcal{D}(\mathfrak{p})$ 
of the prime ideal $\mathfrak{p}:=\mathrm{ker}(\mathrm{res}_p)$ of $\mathrm{O}_\mathcal{R},$
known as the {\em decomposition group of }  $\mathfrak{p}.$
The decomposition group  $\mathcal{D}(\mathfrak{p})$  acts on $O_{\mathcal{R}}/\mathfrak{p}=\F_{\tilde{p}}\cong \F_p^{alg}$ inducing
all the automorphisms of    $\F_p^{alg},$ which determines a canonical surjective homomorphism 
$ \mathcal{D}(\mathfrak{p})\to  \mathcal{G}_{\F_p}$ onto the absolute Galois group of $\F_p$.

\medskip

\subsection{
Question} {\em Is the action of  $\sigma\in  \mathcal{G}_\Q$ on $(\pi^\Spec_{\mathcal{R},1})^{-1}(p)$ over a finite prime continuous?}

\subsection{The fibre of $A_\mathcal{R}$ over the infinite splitting prime }
In this case   $\tilde{p}$ lies over an infinite splitting prime  $\mathbf{p}$ (which is unique up to its first-order type) 
and $\mbox{O}_\mathcal{R}$ contains all integral extensions, we have  that 
$\mbox{O}_\mathcal{R}\subset \Q^{alg}$ is the ring of all integral algebraic numbers,
$\pi_{\mathcal{R}, 1}$ is an isomorphism and
$\F_{\tilde{p}}=\F_\mathbf{p},$ a pseudo-finite field which contains $\Q^{alg}$ and the latter is equal to the subfield  of named elements of  $\F_\mathbf{p}$ (see \ref{compact}).

 Again we will have as many points $\tilde{p}$ over $\mathbf{p}$ as there are  naming homomorphisms $\sigma:\mbox{O}_\mathcal{R}\to \Q^{alg}\subset \F_\mathbf{p}.$ Note that the naming homomorphisms in this case can be identified with automorphisms of the field $\Q^{alg}.$ In other words,  $\mathcal{G}_\Q$ acts freely and transitively on $(\pi^\Spec_{\mathcal{R},1})^{-1}(\mathbf{p})$ 
  $$(\pi^\Spec_{\mathcal{R},1})^{-1}(\mathbf{p})=\tilde{p}\cdot\mathcal{G}_\Q.$$

\subsection{Proposition} {\em The action of a $\sigma\in  \mathcal{G}_\Q$ on the fibre $(\pi^\Spec_{\mathcal{R},1})^{-1}(\mathbf{p})$ over the infinite splitting prime
is continuous.}

{\bf Proof.} The action on each layer $\Spec_K$
is definable by formula (\ref{qq}) which also defines a continuous map according to our definition of topology. $\Box$

 \section{Concluding remarks and further directions}
 
 As was noted in the introduction this version of the structure is the most basic one. We would like to indicate several direction in which the construction and the analysis may develop.
 
 First remark concerns the  similarity with the basic ingredients of Arakelov's geometry. The minimal compact model of the arithmetic plane over $K$ (see \ref{compact}) is quite similar to Arakelov's plane over the projective line.  The limit fibres in our case corresponds to $\Q^{alg},$ whereas in Arakelov's setting it is $\C$ or $\R,$ depending on the number field, and $K$ embeds in the limit fibres  (it is worth noting that in our structure the number of embedding morphisms is $\deg(K:\Q)$ just as Arakelov's theory predicts).
  One way of closing this gap is to consider $\Z$ and the general rings $O_K$ for Galois extension as ``$*$-algebras'', that is with the involution $*,$ complex conjugation  induced by an embedding $K\subset \C$.
   In particular, $\Z$ consists of self-adjoint operators.
 Respectively, we replace in the limit fibres  $\Q^{alg}$ by its completion $\C$ or by its self-adjoint part $\R.$ 
 
 The other connection is with the work of A.Connes and C.Consani \cite{CC}. Our current understanding is that to come to their arithmetic site from our arithmetic plane we need to extend our notion of representations of $\Z$ from $\Z/p,$ for prime $p,$ to the more general $\Z/n,$ for arbitrary $n>1.$ The points of the arithmetic site then can be seen as the limits of such representations.
 
 Finally, the model theorist reader would have already noted that the geometry of our structure is of trivial type in the sense of model theoretic trichotomy (the one by  Connes and Consani is non-trivial  locally modular type). This raises the question if such a geometry can contain any really interesting mathematical information.  To turn this doubt around we note that geometries of trivial  and locally modular types
 may allow much stronger {\em counting functions} than  the non-locally ones. In \cite{Z1984} the first author introduced polynomial invariants of definable subsets in totally categorical theories and in \cite{KMZ} it was proved that the same polynomial invariants in the context of certain combinatorial geometries of trivial type are equivalent to classical graph polynomials of very general kind.     
 
 %\footnote{Standard fibres -- subfields of $\tilde{Q}.$
%Analytic Zariski geometry. Saturated version. Image of $K$ - real and complex. }

\end{document}